\documentclass[11pt,reqno,a4paper]{amsart}
\usepackage{color}
\usepackage{amssymb,amsmath}
\usepackage[latin1]{inputenc}
\usepackage[active]{srcltx}
\usepackage{epsfig,graphics}
\usepackage{graphicx}
\usepackage{exscale,relsize}
\usepackage{textgreek}
\usepackage{psfrag}

\def\N{\mathbb{N}}
\def\R{\mathbb{R}}
\def\m1{{I\!\!M}}

\renewcommand{\to}{\rightarrow}
\newcommand{\pa}{\partial}

\newcommand{\ino}{\int\limits_{\Omega}}

\newcommand{\fo}{\forall}

\newcommand{\rife}[1]{(\ref{#1})}
\newcommand{\ov}[1]{\overline{#1}}
\newcommand{\un}[1]{\underline{#1}}

\newcommand{\sscp}{\scriptscriptstyle}
\newcommand{\dsp}{\displaystyle}

\renewcommand{\dfrac}{\displaystyle\frac}
\newcommand{\finedim}{\hspace{\fill}$\square$}
\newcommand{\intbar}{\mathop{\int\makebox(-15.5,0){\rule[6pt]{.7em}{0.3pt}}%
\kern-6pt}\nolimits}

\newcommand{\ii}{\infty}

\newcommand{\al}{\alpha}

\newcommand{\sg}{\sigma}

\newcommand{\om}{\Omega}
\newcommand{\lm}{\lambda}

\newcommand{\omb}{\ov{\Omega}}

\newcommand{\rl}{\mbox{\Large \textrho}_{\!\sscp \lm}}

\renewcommand{\rho}{\mbox{\Large \textrho}}

\newcommand{\pl}{\psi_{\sscp \lm}}

\newcommand{\zl}{z_{\ssb}}
\newcommand{\zlz}{z_{\sscp \lm,0}}
\newcommand{\wlz}{w_{\sscp \lm,0}}
\newcommand{\el}{E_{\ssb}}

\newcommand{\wl}{w_{\sscp \lm}}

\newcommand{\ul}{u_{\sscp \lm}}

\newcommand{\etb}{\eta_{\sscp \lm}}

\newcommand{\ssb}{\sscp \lm}

\newcommand{\ol}{\ssb,\sscp 0}

\newcommand{\vm}{v_{\sscp \mu}}
\newcommand{\um}{u_{\sscp \mu}}

\newcommand{\mus}{\mu_{\star}}

\newcommand{\prl}{{\textbf{(}\mathbf P_{\mathbf \lm}\textbf{)}}}
\newcommand{\qrl}{{\textbf{(}\mathbf Q_{\mathbf \lm}\textbf{)}}}

\newcommand {\pmu}{\mbox{\rm $(1)_{\mu}$}}

\newtheorem{theorem}{Theorem}[section]
\newtheorem{proposition}[theorem]{Proposition}
\newtheorem{lemma}[theorem]{Lemma}
\newtheorem{corollary}[theorem]{Corollary}
\newtheorem{remark}[theorem]{Remark}
\newtheorem{definition}[theorem]{Definition}
\newcommand{\brm}{\begin{remark}\rm}
\newcommand{\erm}{\end{remark}}
\newcommand{\bdf}{\begin{definition}\rm}
\newcommand{\edf}{\end{definition}}
\newcommand{\bte}{\begin{theorem}}
\newcommand{\ete}{\end{theorem}}
\newcommand{\bpr}{\begin{proposition}}
\newcommand{\epr}{\end{proposition}}
\newcommand{\ble}{\begin{lemma}}
\newcommand{\ele}{\end{lemma}}
\newcommand{\bco}{\begin{corollary}}
\newcommand{\eco}{\end{corollary}}
\newcommand{\beq}{\begin{equation}}
\newcommand{\eeq}{\end{equation}}
\newcommand{\bdm}{\begin{displaymath}}
\newcommand{\edm}{\end{displaymath}}

\newcommand{\graf}[1]{\left\{\begin{array}{ll}#1\end{array}\right.}

\def\sideremark#1{\ifvmode\leavevmode\fi\vadjust{\vbox to0pt{\vss
 \hbox to 0pt{\hskip\hsize\hskip1em \vbox{\hsize2.1cm\tiny\raggedright\pretolerance10000 \noindent #1\hfill}\hss}\vbox to15pt{\vfil}\vss}}}

\begin{document}
\numberwithin{equation}{section}
\parindent=0pt
\hfuzz=2pt
\frenchspacing

\title[]{On the global bifurcation diagram of the Gel'fand problem}

\author[D.B. \& A.J.]{Daniele Bartolucci$^{(1,\ddag)}$, Aleks Jevnikar$^{(2)}$}


\thanks{$^{(1)}$Daniele Bartolucci, Department of Mathematics, University
of Rome {\it "Tor Vergata"}, \\  Via della ricerca scientifica n.1, 00133 Roma,
Italy. E-mail: bartoluc@mat.uniroma2.it}

\thanks{$^{(2)}$ Department of Mathematics, Computer Science and Physics, University of Udine, \\ Via delle Scienze 206, 33100 Udine, Italy. E-mail: aleks.jevnikar@uniud.it}

\thanks{$^{(\dag)}$Research partially supported by:
PRIN project 2015 "{\em Variational methods, with applications to problems in mathematical physics and geometry}",
Consolidate the Foundations project 2015 (sponsored by Univ. of Rome "Tor Vergata") "{\em Nonlinear Differential Problems and their Applications}",
S.E.E.A. project 2018 (sponsored by Univ. of Rome "Tor Vergata"),
MIUR Excellence Department Project awarded to the Department of Mathematics, Univ. of Rome Tor Vergata, CUP E83C18000100006.}

\begin{abstract}
For domains of first kind (\cite{BLin3,CCL}) we describe the qualitative behavior of the global bifurcation diagram of the
unbounded branch of solutions of the Gel'fand problem crossing the origin. At least to our knowledge this is the first result about
the exact monotonicity of the branch of non-minimal solutions which is not just concerned with radial solutions (\cite{Kor}) and/or with symmetric domains (\cite{HK}). Toward our goal we parametrize the branch not by the $L^{\ii}(\om)$-norm of the solutions but
by the energy of the associated mean field problem. The proof relies on a refined spectral analysis of mean field type equations and some surprising properties of the quantities triggering the monotonicity of the Gel'fand parameter.

\end{abstract}
\maketitle

\bigskip

\section{Introduction} \label{sec:intro}
\begin{figure}[h]
\psfrag{L}{$\mu_\star$}
\psfrag{V}{$\mu$}
\psfrag{U}{sgn$(\vm)\|\vm\|_\ii$}
\psfrag{Z}{$(0,0)$}
\includegraphics[totalheight=2in, width=3in]{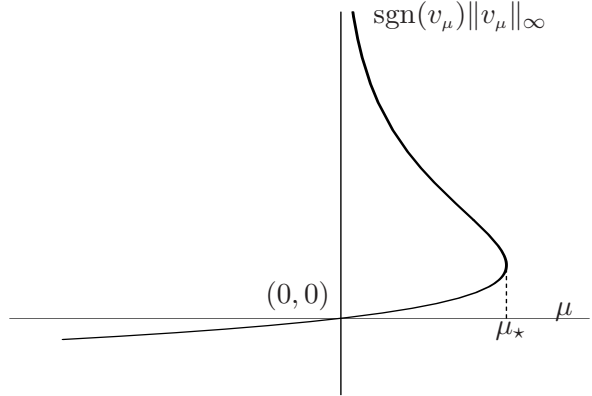}
\caption{The graph of $\Gamma_\ii(B_1)$}\label{fig1}
\end{figure}

\setcounter{equation}{0}
We are concerned with the global bifurcation diagram of solutions of,
$$
\graf{-\Delta v =\mu {\dsp e}^{\dsp v} \qquad   \mbox{in } \om \\ v = 0 \qquad\qquad\quad\!\! \mbox{on } \pa\om}\qquad\qquad \qquad \qquad \pmu
$$
where $\om\subset \R^2$ is any smooth, open and bounded domain and $\mu\in \R$. Problem $\pmu$, also known as the Gel'fand problem (\cite{Gel}), arises in many applications, such as the thermal ignition of gases (\cite{beb}), the dynamics of self-interacting particles (\cite{bav}) and of chemiotaxis aggregation (\cite{suzC}), the statistical mechanics of point vortices (\cite{clmp2}) and of self-gravitating objects with
cylindrical symmetries (\cite{KLB}, \cite{Os}). See also \cite{Bw}, \cite{BKN} and the references quoted therein. A basic question seems unanswered so far concerning the qualitative behavior of the unbounded continuum (\cite{Rab}) of solutions of $\pmu$,
$$
\Gamma_{\ii}(\om)=\Bigr\{(\mu,\vm)\in \R\times C_0^{2,\al}(\,\omb\,):\mbox{\rm $\vm$ solves $\pmu$ for some $\mu\in\R$}\Bigr\},\quad (\mu,\vm)=(0,0)\in \Gamma_{\ii}(\om),
$$
emanating from the origin $(\mu,\vm)=(0,0)$. Under which conditions on $\om$, $\Gamma_{\ii}(\om)$ takes the same form (see Fig. \ref{fig1}) as that  corresponding to a disk $\om=B_R$?  Here $B_R=\{x\in \R^2\,:\,|x|<R\}$ and in this case solutions are radial (\cite{gnn}) and can be evaluated explicitly, see for example
\cite{suz}.

For a general domain, classical results (\cite{CrRab,KK,KC}) show that $\vm$ is a monotonic increasing function of $\mu$ as far as the first eigenvalue of the associated linearized problem is strictly positive. This is the so called branch of minimal solutions which is well understood and naturally described in the
$(\mu,\|\vm\|_{L^{\ii}(\om)})$-plane for $\mu<\mus(\om)\in (0,+\ii)$. In fact, this nice behavior breaks down
at some positive value $\mus(\om)$, which is the least upper bound  of those $\mus$ such that $\pmu$ has solutions for any $\mu<\mus$. In particular, the first eigenvalue of the linearized problem at $\mus(\om)$ is zero and $\mus(\om)$ is known to be a bending point, see \cite{suzC}.\\
The situation for non-minimal solutions on $\Gamma_{\ii}(\om)$, i.e. after the first bending point, is more involved. Besides classical facts (see \cite{Ban1}, \cite{suz} and \cite{Lions} for a complete discussion and references), at least to our knowledge there are only two rather general result
concerning this problem. The first one is in \cite{suz}, where it is shown that, for a certain class of simply connected domains (see Remark \ref{rem1.3} below), $\Gamma_{\ii}(\om)$ is a smooth curve with only one bending point which makes 1-point blow up (\cite{NS90}) as $\mu\to 0^+$. The second one is an unpublished but straightforward corollary of some results in \cite{BLin3}, which implies that the result in \cite{suz}
holds even for a large class of domains with holes, see Remark \ref{rem1.3} below.  Therefore it is natural in this situation to guess that $\Gamma_{\ii}(\om)$ takes the same form shown in Fig. \ref{fig1}. However, a subtle point arises since after the bending, even in this situation where we know that the bifurcation curve cannot bend back to the right, the first eigenvalue of the linearized equation for $\pmu$ is negative (while the second eigenvalue is positive \cite{suz,BLin3}) and then the monotonicity of $\|\vm\|_{\ii}=\|\vm\|_{L^{\ii}(\om)}$, depending in a tricky way on certain changing sign quantities, cannot be taken for granted. Of course, one expects that, as is the case for radial solutions, $\|\vm\|_{\ii}$ is still a monotone function of $\mu$, and well known pointwise estimates (\cite{CLin1,yy}) for blow up solutions suggest that this is the case for $\mu\searrow 0^+$ small enough. At least to our knowledge there are no proofs of this fact. The difficulty is that along non minimal solutions the Morse index is one,
and with only one negative eigenvalue we would like to prove that $\|\vm\|_{\ii}$ is \un{decreasing} in $(0,\mus(\om))$. This is a
kind of inverse maximum principle, since, as mentioned above, if the first eigenvalue is positive then $\|\vm\|_{\ii}$ is increasing.\\ Actually, under some symmetry assumptions on $\om$, by the result in \cite{HK}, for any $m\in (0,+\ii)$ there exists one and only one solution of $\pmu$ such that $\|\vm\|_{\ii}=m$ and $\Gamma_{\ii}(\om)$ is a smooth curve which contains all solutions of $\pmu$. Therefore, for these symmetric domains, the results in \cite{suz} and \cite{HK} together show that indeed $\|\vm\|_{\ii}$ is monotone along $\Gamma_\ii(\om)$, which answer to our question in this case. Finally, it seems that there is no gain in replacing $\|\vm\|_{\ii}$ with other seemingly natural quantities, as for example
$\ino |\nabla \vm|^2$ or either $\ino {\dsp e}^{\dsp \vm}$, as one is always left with the problem of possibly sign changing terms.

\bigskip

We attack this problem here by a new method based on some ideas recently introduced in \cite{Bons}, that is, to parametrize the curve $\Gamma_{\ii}(\om)$ not by $\|\vm\|_{\ii}$ but with the energy $\mathbb{E}(\mu)$ of the associated mean field equation, naturally arising in the Onsager description of two-dimensional turbulence (\cite{clmp2}). See also \cite{Bw} for other results following 
from this approach. Our proof works for domains of "first kind" (see Definition \ref{firstk} below), initially introduced in statistical mechanics (\cite{clmp2}) and then sharpened and fully characterized in \cite{CCL} and in \cite{BLin3}.
For any pair $(\mu,\vm)\in (\R,C_0^{2,\al}(\,\omb\,))$ solving $\pmu$ we define,
$$
\mathbb{E}(\mu)=\graf{\dfrac{1}{2\mu \ino e^{\vm}}\ino\dfrac{ e^{\dsp \vm} }{ \ino e^{\dsp \vm}}\vm, \qquad\; \mu\neq 0,\\ \\
\frac{1}{2|\om|^2}\ino G(x,y)dx dy, \quad \mu=0,}
$$
where $G(x,y)$ is the Green function for $-\Delta$ with Dirichlet boundary conditions.
For later use let us set,

$$
E_0=E_0(\om):=\mathbb{E}(0)=\frac{1}{2|\om|^2}\ino\ino G(x,y)dx dy.
$$
We say that
$f:I\to X$, where $I\subseteq \R$ is an open set and $X$ is a Banach space, is real analytic \cite{but} if for each $t_0\in I$ it admits a power series expansion in $t$, which is totally convergent in the $X$-norm in a suitable neighborhood of $t_0$. Our main result is the following:

\bte\label{thm11} Let $\om$ be a domain of first kind {\rm(}Definition \ref{firstk}{\rm )}. For any $E\in (0,+\ii)$, the equation
$$
\mathbb{E}(\mu)=E \qquad (\,\mu,\vm)\in\Gamma_{\ii}(\om) \qquad\qquad\qquad {\bf (E)}
$$
admits a unique solution $\mu=\mu_{\ii}(E)$. \\ In particular, $\mu_{\ii}:(0,+\ii)\to(-\ii,+\ii)$  and $\left.\vm\right|_{\mu=\mu_{\ii}(E)}:(0,+\ii)\to C^{2,\al}_0(\,\omb\,)$ are real analytic functions of $E$ and
$\left. (\mu,\vm)\right|_{\mu=\mu_{\ii}(E)}$ is a parametrization of $\Gamma_{\ii}(\om)$. Finally
$\mu_{\ii}(E)$ has the following properties:\\
$(i)$ $\mu_{\ii}(E)\to -\ii$ as $E\to 0^+$, $\mu_{\ii}(E_0)=0$, $\mu_{\ii}(E)\to 0^+ $ as $E\to +\ii$;\\
$(ii)$ $\frac{d\mu_{\ii}(E)}{dE}>0$ for $E<E_*$, $\frac{d\mu_{\ii}(E_*)}{dE}=0$, $\frac{d\mu_{\ii}(E)}{dE}<0$ for $E>E_*$, where $E_*=E_*(\om)>E_0(\om)$
is uniquely defined by $E_*(\om)=\mathbb{E}(\mus(\om))$, that is $\mu_{\ii}(E_*)=\mus(\om)$.\\
\ete

\begin{figure}[h]
\psfrag{L}{$\mu_\star$}
\psfrag{V}{$\mu$}
\psfrag{U}{$\mathbb{E}(\mu)$}
\psfrag{E}{$E_*$}
\psfrag{Z}{$(0,E_0)$}
\includegraphics[totalheight=2in, width=3in]{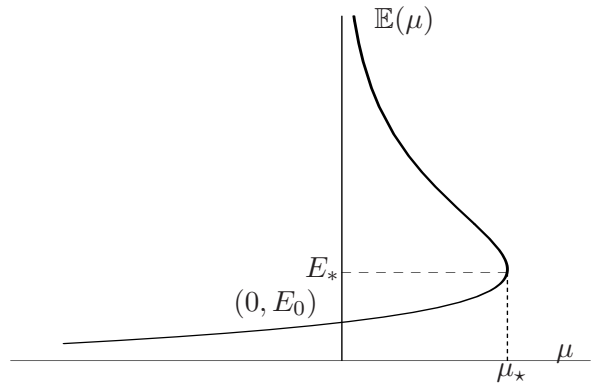}
\caption{The graph of $\Gamma_\ii(\om)$ on domains of first kind}\label{fig2}
\end{figure}

Therefore, on domains of first kind, we have found a global parametrization of $\Gamma_{\ii}(\om)$,
$$
\Gamma_{\ii}(\om)=\left\{(\mu,\vm)\in[0,\mus(\om)]\times C_0^{2,\al}(\,\omb\,)\,:\, \mu=\mu_{\ii}(E),\,E\in (0,+\ii) \right\},
$$
which takes the form depicted in Fig. \ref{fig2}, as claimed. At least to our knowledge, this is the first global result (i.e. including non-minimal solutions) about the monotonicity of the bifurcation diagram for an elliptic equation with superlinear growth in dimension $n= 2$, which is not just concerned with radially symmetric solutions (\cite{Kor}, \cite{KLO}),  and/or with domains sharing some kind of symmetries (\cite{HK}). The situation in higher dimension is far more subtle, see for example \cite{JL} and more recently \cite{Dan2,Dan3,DanF,Kor2}.
\bigskip

Although its definition in the context of the Gel'fand problem looks rather unnatural, it turns out that indeed $\mathbb{E}(\mu)$ is just the energy in the Onsager mean field model, when expressed as a function of $\mu$.  For $\lm\in (-\ii,8\pi)$ we consider the mean field equation (\cite{clmp2}),
$$
\graf{-\Delta \pl =\dfrac{e^{\lm\pl}}{ \ino e^{\lm\pl}} \quad \om \\ \pl = 0 \qquad \qquad\pa\om}\qquad\qquad \qquad \qquad \prl
$$
To simplify the notations, for fixed $\lm\in (-\ii,8\pi)$, here and in the rest of this paper we set,
$$
\rl=\dfrac{e^{\lm\pl}}{ \ino e^{\lm\pl}}, \qquad <f>_{\ssb}=\ino \rl f, \qquad f_0=f-<f>_{\ssb}.
$$
The energy associated with the density $\rl$ is by definition (see \cite{clmp2}),
$$
\mathcal{E}(\rl)=\frac12\ino\ino  \rl(y)G(x,y)\rl(x)dydx.
$$
Clearly, if $\pl$ solves $\prl$, then $\vm=\lm\pl$ is a solution of $\pmu$, for some $\mu\in \R$ which satisfies,
$$
\mu=\mu_{\ssb}=\frac{\lm}{\ino e^{\lm\pl}},
$$
and then from $\prl$ we see that,
$$
\el:=\mathbb{E}(\mu_{\ssb})=\mathcal{E}(\rl)=\frac{1}{2} <\pl>_{\ssb}=\frac{1}{2}\ino |\nabla \pl|^2.
$$
Therefore, in particular the energy $\mathbb{E}(\mu_{\ssb})$ is just the Dirichlet energy of $\pl$ when expressed in terms of $\mu$.\\
The uniqueness of solutions of $\prl$ is easy to prove for $\lm< 0$, see Proposition \ref{positive} below. On the other side it is well known that $\prl$ admits a unique solution for any $\lm\in (0,8\pi)$ (see \cite{suz} for simply connected domains and \cite{BLin3} for general domains).
The existence/non existence problem for $\lm=8\pi$ is a subtle issue, since $\prl$ is critical with respect to the Moser-Trudinger inequality \cite{moser}. For a complete discussion of this problem see \cite{CCL} for simply connected domains and \cite{BLin3} for general domains. This is why we need
the following definition,

\bdf\label{firstk} {\it A domain $\om$ is of first kind if {\rm $\prl$} has no solution for $\lm=8\pi$ and is of second kind otherwise.}
\edf

\brm\label{rem1.3} {\it Concerning the discussion in the first part of the introduction, we point out that 
it has been proved in \cite{CCL} and \cite{BLin3} that
$\om$ is of first kind if and only if  the unique solutions of {\rm $\prl$} for $\lm<8\pi$ blow up {\rm (}\cite{bm,NS90}{\rm )}, as $\lm\to 8\pi^-$. As a consequence, it turns out that the domains considered in Theorem 2 in \cite{suz} are exactly the simply connected domains of first kind. However the spectral estimates in \cite{suz} has been extended {\rm (}and improved up to $\lm=8\pi${\rm )} in \cite{BLin3} to the case of domains with non trivial topology. As an immediate consequence,
the results in \cite{suz} mentioned in the first part of the introduction hold for any domain of first kind.\\ It is well known that any disk $\om=B_R$ is of first kind. Actually any regular polygon is of first kind, see \cite{CCL}.
It has been proved in \cite{BdM2} that there exists a universal constant $I>4\pi$ such that any convex domain whose isoperimetric
ratio is larger than $I$ is of second kind. If $\om_{a,b}$ is a rectangle of sides $a\leq b$, then there exists $\xi\in(0,1)$ such that
$\om_{a,b}$ is of first kind if and only if $\frac{a}{b}\leq \xi$, see \cite{CCL}. If $\om_r=B_1\setminus \ov{B_r(x_0)}$ with $x_0\in B_1$,
$x_0\neq 0$, then there exists $r_0<\min\{|x_0|, 1-|x_0|\}$ such that $\om_r$ is of first kind for any $r<r_0$, see \cite{BLin3}.
We refer to \cite{CCL} and \cite{BLin3}  for other equivalent characterizations of domains of first kind and a complete discussion concerning this point. Among other things it is proved there that the set of domains of first kind is closed in the $C^1$-topology and that domains of first kind need not be symmetric. More recently it has been proved in \cite{BMal} that the set of domains of first kind has non-empty interior in the $C^{1}$-topology and that the set of simply connected $C^{1}$-domains of first kind is contractible.}
\erm

\brm{\it
 It is well known {\rm(}\cite{CLin2},\cite{KMdP}{\rm)} that if $\om$ is not simply connected, then there exist countably many distinct families of blow up solutions of $\pmu$ as $\mu\to 0^+$. These families of blow up solutions are unique {\rm(}\cite{BJLY3}{\rm)} and nondegenerate {\rm(}\cite{GOS}{\rm)} under suitable nondegeneracy assumptions. Therefore, it is not true in general that $\Gamma_\ii(\om)$ contains all solution of $\pmu$. However, as remarked above,
 this is the case for domains with certain symmetries, see \cite{HK}.}
\erm

\bigskip

 Let us sketch the argument of the proof. We describe $\Gamma_\ii(\om)$ by using the unique solutions of $\prl$ for $\lm\in (-\ii,8\pi)$. Any solution $\pl$ of $\prl$ yield a solution of $\pmu$ with $\mu_{\ssb}=\lm(\ino e^{\lm\pl})^{-1}$. This might seem not a good point of view since, $\prl$ being a constrained problem, the associated linearized operator $L_{\ssb}$ (see \rife{2.1} and \rife{lineq0} below) is more difficult to analyze. In general the first eigenvalue ${\sg}_{1,\lm}$ is not simple, the first eigenfunctions may change sign and the positivity of $\sg_{1,\lm}>0$ does not imply in general that the maximum principle holds, see \cite{B2} or the Appendix below for an example of this sort. Therefore, even if we know that ${\sg}_{1,\lm}>0$ (see \cite{CCL} and \cite{BLin3}), it is not clear how to use this information to establish the monotonicity of $\mu_{\ssb}$. However, as recently observed in \cite{Bons}, it is possible to set up the spectral theory relative to the linearization of $\prl$ and build a complete set of eigenfunctions which span the space of functions of vanishing mean.
  All the main steps of the proof rely on this carefully defined spectral setting. The first eigenvalue $\sg_{1,\lm}$ is strictly positive for $\lm\in (-\ii,8\pi)$ and the set of solutions of $\prl$ can be shown to be locally a real analytic curve $(\lm,\pl)$ with no bifurcation points. A crucial fact which follows from $\sg_{1,\lm}>0$ is that the energy $\el$  is a (real analytic) strictly increasing function of $\lm$. This follows from a careful analysis of the Fourier modes of $L_{\ssb}$. Therefore, to understand the monotonicity of $\mu$ as a function of $E$, it is enough to evaluate the sign of $\frac{d \mu_{\ssb}}{d\lm}$. The hard part is to show that there exists $\lm_*\in (0,8\pi)$ such that $\frac{d \mu_{\ssb}}{d\lm}>0 \iff \lm<\lm_*$. Indeed, a major problem arises in the proof of $\frac{d\mu_{\ssb}}{d\lm}<0$ along the non-minimal branch of solutions, that is for $\lm>\lm_*$. We solve this problem by two non-trivial facts about the quantity which controls the sign of $\frac{d \mu_{\ssb}}{d\lm}$, which is $g(\lm)$ in Lemma \ref{lem1.3} below. First of all, still by exploiting $\sg_{1,\lm}>0$ and based on the sign of $g(\lm)$, we obtain an integral version of the maximum principle for a particular test function.
  This is not at all obvious since, as mentioned above, $\sg_{1,\lm}>0$ does not imply that the maximum principle holds for $L_{\ssb}$. The second fact is a remarkable formula for $g(\lm)$: it turns out that it satisfies a first order linear non-homogeneous O.D.E. (see \rife{cruxg} below). These two properties together yield the desired form of the inverse maximum principle, as
  we conclude in this way that $\frac{d \mu_{\ssb}}{d\lm}$ changes sign only once in $(0,8\pi)$.\\

\bigskip

The description of $\Gamma_{\ii}(\om)$ on domains of second kind is more difficult. Indeed, solutions on a certain part of $\Gamma_{\ii}(\om)$ correspond to solutions of $\prl$ with $\lm>8\pi$, a region where solutions are not unique and $\sg_{1,\lm}$ 
in general is not positive. Therefore
it is not easy to understand the monotonicity of $\el$, see \cite{Bons} for some results concerning this point.

\bigskip

\bigskip

This paper is organized as follows. In section \ref{sec:spectral} we first introduce the spectral setting for $L_{\ssb}$ and then collect some important preliminary results concerning the linearized mean field equation \rife{lineq0}. In section \ref{sec:mono} as a first step toward the proof of the main result we deduce the monotonicity of the energy. Then, in section \ref{sec:proof.thm} we prove the main Theorem \ref{thm11} postponing the proof of the key Lemma \ref{lem1.3} to section \ref{sec:proof.lem}. Finally, further discussions and an explicit example about the operator $L_{\ssb}$ are given in the Appendix.

\bigskip
\bigskip
\section{Spectral decomposition of linearized mean field type equations} \label{sec:spectral}

For any $\pl$ solving $\prl$, we introduce the linearized operator,
\beq\label{2.1}
L_{\ssb}\phi:=-\Delta \phi-\lm\rl \phi_0,\quad \phi\in H^1_0(\om)
\eeq
where we recall that
$$
\phi_0=\phi \;- <\phi>_{\ssb}.
$$

We say that $\sg=\sg(\lm,\pl)\in\R$ is an eigenvalue of the linearized operator \rife{2.1} if the equation,
\beq\label{lineq0}
-\Delta \phi-\lm \rl \phi_0=\sg\rl \phi_0,
\eeq
admits a non-trivial weak solution $\phi\in H^1_0(\om)$. This definition of the eigenvalues requires some comments.
Let $\pl$ be a fixed solution of $\prl$ and let us define,
$$
Y_0:=\left\{ \varphi \in \{L^2(\om),<\cdot,\cdot >_{\ssb}\}\,:\,\ino\rl \varphi=0\right\},
$$
where $<\cdot,\cdot>_{\ssb}$ denotes the scalar product $<f,g>_{\ssb}:=<fg>_{\ssb}$.
Let us also define
$$
T(\phi):=G[\rl \phi], \quad \phi\in L^2(\om),\mbox{ where }G[f]=\ino G(x,y)f(y)dy.
$$
By the results in \cite{bm,yy} and standard elliptic regularity theory we see that $\rl$ is a smooth function for $\lm<8\pi$, so these definitions are well posed. Clearly $Y_0$ is an Hilbert space, and,
since $T(Y_0)\subset W^{2,2}(\om)$, then it is not difficult to see that the linear operator,
\beq\label{T0}
T_0:Y_0\to Y_0,\;T_0(\varphi)=G[\rl\varphi]-<G[\rl\varphi]>_{\ssb},
\eeq
is self-adjoint and compact. As a consequence,
standard results concerning the spectral decomposition of self-adjoint, compact operators on Hilbert spaces show that
$Y_0$ is the Hilbertian direct sum of the eigenfunctions of $T_0$, which can be represented as
$$
\varphi_k=\phi_{k,0}:=\phi_k-<\phi_k>_{\ssb},\;k\in\N=\{1,2,\cdots\},
$$
$$
Y_0=\overline{\mbox{Span}\left\{\phi_{k,0},\;k\in \N\right\}},
$$
for some $\phi_k\in H^1_0(\om)$, $k\in\N=\{1,2,\cdots\}$.
In fact $\varphi_k$ is an eigenfunction whose eigenvalue is $\mu_k=\frac{1}{\lm+\sg_k}\in \R\setminus\{0\}$,
that is,
$$
\varphi_k=(\lm+\sg_k) \left(G[\rl\varphi_k]-< G[\rl\varphi_k]>_{\ssb}\right),
$$
if and only if the function $\phi_k$,
$$
\phi_k:=(\lm+\sg_k) G[\rl\varphi_k],
$$
is in $H^1_0(\om)$ and weakly solves,
\beq\label{lineq}
-\Delta \phi_k= (\lm+\sg_k) \rl \phi_{k,0}\quad\mbox{ in }\quad \om.
\eeq

At this point, standard arguments in the calculus of variations show that,
\beq\label{4.1}
\sg_1=\sg_1(\lm,\pl)=\inf\limits_{\phi \in H^1_0(\om)\setminus \{0\}}
\dfrac{\ino |\nabla \phi|^2 - \lm \ino \rl\phi_0^2 }{\ino  \rl\phi_0^2}.
\eeq
The ratio in the right hand side of  \rife{4.1} is well defined in $H^1_0(\om)$ by the Jensen inequality, which implies that
$$
\ino  \rl\phi_0^2=<\phi^2>_{\ssb}-<\phi>^2_{\ssb}\geq 0,
$$
where the equality holds if and only if $\phi\equiv 0$ a.e. in $\om$. Higher eigenvalues are defined inductively via the variational problems,
\beq\label{4.1.k}
\sg_k=\sg_k(\lm,\pl)=\inf\limits_{\phi \in H^1_0(\om)\setminus \{0\},<\phi_0,\phi_{k,0}>_{\ssb}=0,\; m\in \{1,\ldots, k-1\}}
\dfrac{\ino |\nabla \phi|^2 - \lm \ino \rl\phi_0^2 }{\ino  \rl\phi_0^2}.
\eeq
Obviously a base  $\{\phi_{i,0}\}_{i\in\N}$ of $Y_0$ can be constructed to satisfy,
\beq\label{orto}
<\phi_{i,0}, \phi_{j,0}>_{\ssb}=0,\;i\neq j.
\eeq
The eigenvalues form a countable nondecreasing
sequence $\sg_1(\lm,\pl)\leq \sg_2(\lm,\pl)\leq....\leq \sg_k(\lm,\pl)\leq...\,$, where in particular,
\beq\label{lineq2}
\lm+\sg_k\geq \lm +\sg_1>0,\;\fo\,k\in \N,
\eeq
the last inequality being an immediate consequence of \rife{4.1}.\\
Obviously, by the Fredholm alternative,
\beq\label{invert}
\mbox{\rm if $0\notin \{\sg_j\}_{j\in \N}$, then $I-\lm T_0$ is an isomorphism of $Y_0$ onto itself}.
\eeq
Finally, any $f \in L^{2}(\om)$ admits the Fourier series expansion,
\beq\label{Four}
f=\al_{0}+\sum\limits_{j=1}^{+\ii} \al_{j}\phi_{j,0},\quad \ino \rl \phi^2_{j,0}=1,
\eeq
where
$$
\al_{j}=\al_{j}(f)=<\phi_{j,0},f>_{\ssb}.
$$

For later use we define $\widehat{\sg}_1$ to be the standard first eigenvalue defined by,
\beq\label{sghat}
\widehat{\sg}_1=\inf\limits_{\phi \in H^1_0(\om)\setminus \{0\}}\dfrac{\ino |\nabla \phi|^2 - \lm \ino \rl(\phi^2-<\phi>_{\ssb}^2) }{\ino  \rl\phi^2}.
\eeq

A fundamental spectral estimate will be used which follows from the results in \cite{suz}, as later improved in \cite{CCL} and \cite{BLin3}:\\

{\bf Theorem A.}(\cite{BLin3}, \cite{CCL}, \cite{suz})
{\it Let $\om$ be a smooth and bounded domain. For any $\lm\in [0,8\pi)$, the first eigenvalue of
\eqref{2.1} is strictly positive, that is,  $\sg_1(\lm,\pl)>0$, for each $\lm\in[0,8\pi)$.}

\bigskip

The curve,
$$
(-\ii,8\pi) \supseteq I \ni \lm \mapsto (\lm,\pl)\in (-\ii,8\pi)\times C_0^{2,\al}(\ov{\om}),
$$
from an open interval $I\subseteq (-\ii,8\pi)$ is said to be analytic if for each $\lm_0\in I$, $\pl$ admits a power series expansion as a function of $\lm$, which is totally convergent in the $C_0^{2,\al}(\ov{\om})$-norm in a suitable neighborhood of $\lm_0$.
We have the following,
\bpr\label{positive} Let $\om$ be a smooth and bounded domain. \\
$(i)$ For any $\lm\in (-\ii,8\pi)$ it holds $\sg_1=\sg_1(\lm,\pl)>0$.\\
$(ii)$ If $\pl$ solves \mbox{\rm $\prl$} for some $\lm=\lm_0\in (-\ii,8\pi)$, then  there exists an open neighborhood of $(\lm_0,\psi_{\sscp \lm_0})$ in the product topology, $I\times B=\mathcal{U}\subset
\R\times C^{2,\al}_0(\ov{\om})$,  such that the set of solutions of
{\rm $\prl$} in $\mathcal{U}$ is an analytic curve $I\ni\lm\mapsto (\lm,\pl)\in \mathcal{U}$, for
suitable neighborhoods $I$ of $\lm_0$ in $\R$ and $B$ of $\psi_{\sscp \lm_0}$ in $C^{2,\al}_0(\ov{\om})$.\\
$(iii)$ For any $\lm\in (-\ii,8\pi)$ there exists a unique solution of \mbox{\rm $\prl$}.
\epr
\proof
If $\pl$ solves $\prl$ for some $\lm\in (-\ii,8\pi)$, then Theorem A and \rife{lineq2} immediately imply that $\sg_1=\sg_1(\lm,\pl)>0$, which proves $(i)$. As a consequence of $(i)$, \rife{invert} and the analytic implicit function theorem (\cite{bdt, but}), it
can be shown by standard arguments that $(ii)$ holds, see Lemma 2.4 in \cite{Bons}. Therefore also $(ii)$ is proved.\\
Putting $\ul=\lm\pl$, then solutions of $\prl$ are critical points of
$$
J_\lm(u)=\frac12 \ino |\nabla u|^2-\lm \log\left(\ino e^{u}\right), \;u\in H^1_0(\om).
$$
The existence of at least one minimizer for $J_\lm$ for any $\lm\in (-\ii,8\pi)$ is a well known consequence of the Moser-Trudinger inequality (\cite{moser}).
The uniqueness of solutions of $\prl$ for $\lm=0$ is trivial, while for $\lm \in (0,8\pi)$ it has been proved in \cite{BLin3}, \cite{CCL}, \cite{suz}.
For $\lm<0$ the second variation of $J_\lm$ is strictly positive definite. Indeed, by \rife{sghat}, the first eigenvalue of the associated quadratic form, which is $\widehat{\sg}_{1,\lm}$, is strictly positive. Therefore $J_\lm$ admits at most one critical point, which concludes the proof of $(iii)$.
\finedim

\bigskip

\section{Monotonicity of the energy} \label{sec:mono}
The proof of Theorem \ref{thm11} relies on Proposition \ref{positive}, the following result about the global structure of the bifurcation diagram of $\prl$ for $\lm\in(-\ii,8\pi)$ and on the asymptotic behavior of the energy. This is an extension of some ideas first introduced in \cite{Bons}.
\bpr\label{pr1.2}
Let $\om$ be a domain of first kind.
For each $\lm\in (-\ii,8\pi)$ there exists a unique solution $\pl$ of \mbox{\rm $\prl$} and
$$
\mathcal{G}_{8\pi}=\left\{(\lm,\pl)\,:\, \lm\in (-\ii,8\pi),\,\pl\in C_0^{2,\al}(\ov{\om}) \mbox{ solves }\mbox{\rm $\prl$} \right\},
$$
is a real analytic curve. In particular
\beq\label{elm}
\el=\mathbb{E}(\mu_{\ssb})=\frac{1}{2}\ino \mbox{\rm $\rl$}\pl=\frac{1}{2} <\pl>_{\ssb}=\ino |\nabla \pl|^2,
\eeq
is real analytic in $(-\ii,8\pi)$ and satisfies
\beq\label{elm1}
E_{\ssb}\to 0^+ \mbox{ as } \lm\to -\ii, \, E_{\ssb}\to +\ii, \mbox{ as } \lm\to 8\pi^-,\,\left.E_{\lm}\right|_{\lm=0}=E_{ 0},
\eeq
and
\beq\label{delm}
\frac{d E_{\ssb}}{d\lm}>0,\forall\,\lm\in(-\ii,8\pi).
\eeq
\epr
\proof The claim about the uniqueness and regularity of $\pl$ is just Proposition \ref{positive}$(ii)-(iii)$. It is easy to see from $\prl$ that
\rife{elm} holds. Since $\pl$ is real analytic in $(-\ii,8\pi)$, then $\el$, being the composition of real analytic functions, is also real analytic as a function
$\lm$ in $(-\ii,8\pi)$, see for example Theorem 4.5.7 in \cite{but}.\\
Next, let us prove the inequality in \rife{delm}. Since $\pl$ is real analytic, then in particular it is differentiable as a function of $\lm$, and
then by  $\prl$ we conclude that $\etb=\frac{d\pl}{d\lm}\in C^{2,\alpha}_0(\,\omb\,)$ is a classical solution of,
\beq\label{1b1}
-\Delta \etb =\rl\psi_{\ol}+\lm \rl \eta_{\ol}\quad \mbox{ in }\quad\om
\eeq
where
$$
\psi_{\ol}=\pl -<\pl>_{\ssb}\mbox{ and }\eta_{\ol}=\etb -<\etb>_{\ssb}.
$$
By using \rife{elm}, we also conclude that,

\beq\label{8.12.2}
\frac{d \el}{d \lm} =\ino (\nabla \etb,\nabla \pl)=-\ino \etb (\Delta \pl)=<\etb>_{\ssb}.
\eeq
Next, by using $\prl$ and \eqref{1b1}, we have,
\beq\label{2b1}
<\etb>_{\ssb}=\ino\rl\etb =\ino -(\Delta \pl)\etb=\ino -\pl(\Delta \etb)=<\psi_{\ol}^2>_{\ssb}+\lm<\psi_{\ol} \eta_{\ol}>_{\ssb}.
\eeq

Let,
$$
\psi_{\ol}=\sum\limits_{j=1}^{+\ii}\al_j\phi_{j,0}, \quad
\eta_{\ol}=\sum\limits_{j=1}^{+\ii}\beta_j\phi_{j,0},
$$
be the Fourier expansions \rife{Four} of $\psi_{\ol}$ and $\eta_{\ol}$.
Multiplying \rife{1b1} by $\phi_{j,0}$, integrating by parts and using \rife{lineq0}, we conclude that,
\beq\label{lamq31}
\sg_{j}\ino \rl\phi_{j,0}\eta_{\ol} =
\ino\rl\phi_{j,0}\psi_{\ol},\mbox{ that is }\sg_j\beta_j=\al_j,
\eeq
where $\sg_{j}=\sg_{j}(\lm,\pl)$. As a consequence, in view of  Proposition \ref{positive}$(i)$, for any $\lm\in (-\ii,8\pi)$ we conclude that,
\beq\label{event1}
<\psi_{\ol},\eta_{\ol}>_{\lm}=\sum\limits_{j=1}^{+\ii} \al_j\beta_j=
\sum\limits_{j=1}^{+\ii} \sg_{j}(\beta_j)^2\geq \sg_{1}\sum\limits_{j=1}^{+\ii}(\beta_j)^2\geq \sg_{1}<\eta_{\ol}^2>_{\lm},
\eeq
and
\beq\label{event2}
<\psi_{\ol},\eta_{\ol}>_{\lm}=\sum\limits_{j=1}^{+\ii} \al_j\beta_j=
\sum\limits_{j=1}^{+\ii} \frac{\al_j^2}{\sg_{j}}\leq \frac{1}{\sg_{1}}\sum\limits_{j=1}^{+\ii}{\al_j^2}= \frac{1}{\sg_{1}}<\psi_{\ol}^2>_{\lm},
\eeq
By using \rife{8.12.2}, \rife{2b1}, \rife{event1}, \rife{event2} and Proposition \ref{positive}$(i)$, we finally obtain,
\beq\label{deng1}
 \frac{d \el}{d \lm}\geq <\psi_{\ol}^2>_{\lm}+{\lm}{\sg_{1}}<\eta_{\ol}^2>_{\lm}>0,\,\forall \lm\in [0,8\pi),
\eeq
and
$$
\frac{d \el}{d \lm}\geq <\psi_{\ol}^2>_{\lm}+\frac{\lm}{\sg_{1}}<\psi_{\ol}^2>_{\lm}=\frac{\sg_1+\lm}{\sg_1}<\psi_{\ol}^2>_{\lm}>0,\,\forall \lm\in (-\ii,0),
$$
where in the last inequality we used \rife{lineq2}.
The last two inequalities prove the inequality in \rife{delm}.\\
Next let us prove that $\el \to +\ii$ as $\lm\to 8\pi^{-}$. Indeed, if by contradiction there exist
$\lm_n\to 8\pi^{-}$ and $\psi_n=\psi_{\sscp \lm_n}$, such that $E_{\lm_n}\leq C$, $\forall\,n\in\N$, then by \rife{elm} we would conclude that $\ino |\nabla \psi_n|^2\leq C$ for any $n\in\N$. Therefore, since $\pl$ is positive for $\lm>0$, by the Jensen and Moser-Trudinger inequalities (\cite{moser}), we would also conclude that $\rl$ with $\lm=\lm_n$ is uniformly bounded in $L^p(\om)$ for any $p>1$. At this point, by standard elliptic estimates, we could obtain a uniform bound for $\psi_n$ as $\lm_n\to 8\pi^{-}$ and then in particular, passing to the limit along a subsequence in $\prl$, we could conclude that
$\prl$ admits a solution for $\lm=8\pi$. This is impossible since $\om$ is of first kind and so by Definition \ref{firstk} there is no solution of $\prl$ for $\lm=8\pi$. This contradiction proves that $\el\to +\ii$ as $\lm\to 8\pi^{-}$. Obviously, by the uniqueness of solutions,
we also have $\left.E_{\lm}\right|_{\lm=0}=E_{ 0}$.\\
We are just left with the proof of $\el\to 0^+\mbox{ as }\lm\to -\ii$. By contradiction, in view of the monotonicity of $\el$, we can assume that $\el\to E_{\ii}>0\mbox{ as }\lm\to -\ii$. By using once more the monotonicity of $\el$ and $\el\to +\ii$ as $\lm\to 8\pi^{-}$, we conclude that there is no solution of $\prl$ with $\lm\in (-\ii,8\pi)$ whose energy is less than $E_{\ii}$. This is impossible since a well known result in \cite{clmp2} states that for any $E>0$ there is a solution  of $\prl$ with $\lm<8\pi$ whose energy is $E$.
\finedim

\bigskip

\section{The proof of  Theorem \ref{thm11}} \label{sec:proof.thm}

\emph{The proof of  Theorem \ref{thm11}}.\\
The crux of the proof of Theorem \ref{thm11} is to control the monotonicity of $\mu$ as a function of $E$. The first part of this
 task has been accomplished in Proposition \ref{pr1.2}, which is based on the study of solutions of $\prl$. To conclude the proof we will analyze the monotonicity of $\mu$ as a function of $\lm$. Actually, it is easier to do this in terms of $\ul=\lm\pl$ which solves
$$
\graf{-\Delta \ul =\lm\dfrac{e^{\ul}}{ \ino e^{\ul}} \quad \om \\ \ul = 0 \qquad \qquad\pa\om}\qquad\qquad \qquad \qquad \qrl
$$
if and only if $\pl$ solves $\prl$. Clearly, by Proposition \ref{pr1.2}, $\ul$ is a real analytic function of $\lm$, which in particular shows that
this correspondence holds also for $\lm=0$, in the sense that $\ul\equiv 0$ if and only if $\lm=0$
and $\pl$ is the unique solution of $\prl$ with $\lm=0$.

\bigskip
Since, by Proposition \ref{positive}$(iii)$, for
 $\lm\in (-\ii,8\pi)$ $\qrl$ admits a unique solution, then we can define the map $\lm \mapsto \mu_{\ssb}$ as follows,
\beq\label{mul}
\mu_{\ssb}=\dfrac{\lm}{\ino {\dsp e}^{\dsp \ul}}.
\eeq
Clearly $\mu_{\ssb}$ is well defined as a function of $\lm$ for $\lm\in (-\ii,8\pi)$ and obviously, for fixed $\lm$,
the pair $(\mu, \vm)$, where $\mu=\mu_{\ssb}$ and $\vm=\ul$  solves $\pmu$. Although not needed here, it is understood that $\left.(\mu,\vm)\right|_{\mu_{\ssb}}$, for $\lm>0$ small enough, is just a parametrization of a portion of the minimal branch 
(\cite{CrRab},\cite{KK}) of solutions of $\pmu$ for $\mu>0$ small enough.\\
Since $\ul$ is real analytic in $(-\ii,8\pi)$, then, in view of \rife{mul},
we also have that $\mu_{\ssb}$ is a real analytic function of $\lm$ and it holds,
\beq\label{050717.1}
\frac{d \mu_{\ssb}}{d\lm}={\left(\ino e^{\dsp \ul}\right)^{-2}}{\left(\ino e^{\dsp \ul}-\lm\ino  e^{\dsp \ul} \zl\right)},\quad \lm\in (-\ii,8\pi),
\eeq
where
$$
\zl=\frac{d \ul}{d\lm}=\ul^{'}.
$$
Therefore we readily find that,

$$
\left(\ino e^{\dsp \ul}\right)\frac{d \mu_{\ssb}}{d\lm}=1-\lm <\zl>_{\ssb},\quad \lm\in (-\ii,8\pi).
$$

\bigskip

The crux of the proof of Theorem \ref{thm11} is the following,

\ble\label{lem1.3} Let
$g(\lm):=1-\lm <\zl>_{\ssb}$, $\lm \in (-\ii,8\pi)$. There exists $\lm_*\in [4\pi,8\pi)$ such that
$g(\lm)>0$ for $\lm\in (-\ii,\lm_*)$, $g(\lm_*)=0$, $g(\lm)<0$ for $\lm\in (\lm_*,8\pi)$.
In particular $g(0)=1$ and $g(\lm)\to -\ii$ as $\lm\to 8\pi^{-}$.
\ele

\bigskip
\bigskip

We first conclude the proof of Theorem \ref{thm11} and then get back in the next section to that of Lemma \ref{lem1.3}.\\

Let $E_{\ssb}$ be as defined in Proposition \ref{pr1.2} and
$E_*=\left.E_{\ssb}\right|_{\lm=\lm_*}$. By \rife{delm} in Proposition \ref{pr1.2}, we see that $E_{\ssb}$ is a strictly increasing and analytic function of $\lm$ in $(-\ii,8\pi)$ and in particular that $E=E_*$ if and only if $\lm=\lm_*$. Therefore, it is well defined the inverse of $\el$, $\lm_\ii=\lm_\ii(E)$,
and we define
$$
\mu_{\ii}(E)=\left.\mu_{\ssb}\right|_{\lm=\lm_\ii(E)}.
$$

In particular, since $\mu_{\ssb}$ is a real analytic function of $\lm$, and since $\el^{'}$ is strictly positive, then $\lm_\ii$ is a real analytic function of $E$ and so does $\mu_{\ii}$.\\
As a consequence of Lemma \ref{lem1.3} we conclude that
$\frac{d \mu_{\ssb}}{d\lm}>0$ iff $\lm<\lm_*$, where $\lm_*$ is the
unique value of $\lm\in [4\pi,8\pi)$ such that $1-\lm <\zl>_{\ssb}=0$.
Next, by \rife{elm1} in Proposition \ref{pr1.2}, we have $E_{\ssb}\to 0^+$ as $\lm\to -\ii$. Since solutions of $\pmu$ are uniformly bounded in
$[\widehat{\mu},0]$ for any $\widehat{\mu}<0$, since $\pmu$ admits a unique solution for any $\mu<0$, and by using the definition
of $\mathbb{E}(\mu)$, then it is not difficult to see that necessarily $\mu_{\ii}\to -\ii$ as $E\to 0^+$. Also, since $\mu_{\ssb}\to 0$ as $\lm\to 0$, and
$$
\mathbb{E}(\mu_{\ssb})=\el\to E_0, \mbox{ as } \lm \to 0,
$$
then $\mu_\ii(E_0)=0$. By \rife{elm1} in Proposition \ref{pr1.2} we also have $\mathbb{E}(\mu_{\ssb})=E_{\ssb}\to +\ii$ as $\lm\to 8\pi^-$ and by Lemma \ref{lem1.3} $\frac{d \mu_{\ssb}}{d\lm}<0$ iff $\lm>\lm_*$. Therefore $\mu_\ii$ decreases monotonically in $(E_*,+\ii)$ and then we also conclude that $\mu_\ii(E)\to 0^+$ as $E\to +\ii$. Indeed, this is an immediate consequence of \rife{mul} and of the following well known estimate for blow up solutions of mean field type equations, $\ino e^{\ul}\to +\ii$ as $\lm\to 8\pi^{-}$, see for example \cite{yy}.

As a consequence, to conclude the proof, it just remains to show that $E_*=E_*(\om)=\mathbb{E}(\mus(\om))$, that is $\mu_{\ii}(E_*)=\mus(\om)$. Clearly, by the definition of $\mus(\om)$, we have
$\mu_{\ii}(E_*)\leq \mus(\om)$. Therefore we are left to prove the following:\\
{\bf Claim: $\mu_{\ii}(E_*)\geq \mus(\om)$.}\\
Let us assume by contradiction that $\mu_{\ii}(E_*)<\mus(\om)$.
As $E$ increases above $E_*$, then, because of \rife{delm}, we see that $\lm$ is strictly increasing. Therefore $\lm>\lm_*$ for any $E>E_*$ and in particular
$\mu_{\ssb}<\mu_{\ii}(E_*)<\mus(\om)$ for any $\lm\in (0,8\pi)$. As a consequence there exists some open interval $I\subset (0,8\pi)$ such that
$\{(\mu_{\ssb},v_{\mu_{\ssb}})\}_{\lm\in I}\cap
\Gamma_{\ii}(\om)= \emptyset$.

Therefore the solutions $(\mu_{\ssb},v_{\mu_{\ssb}})$ form another smooth branch with no bifurcation points, say  $\mathcal{M}_{8\pi}$,
which however obviously emanates from $\left.(\mu_{\ssb},v_{\mu_{\ssb}})\right|_{\ssb=0}=(0,0)$.

However $\Gamma_{\ii}(\om)$ is a smooth branch with no bifurcation points as well which emanates from $(\mu,\um)=(0,0)$, and this is
a contradiction since then $(0,0)$ should be a bifurcation point where the two branches $\mathcal{M}_{8\pi}$ and $\Gamma_{\ii}(\om)$ meet.
\finedim
\bigskip

\section{The proof of Lemma \ref{lem1.3}} \label{sec:proof.lem}
To simplify the notations, in this section we set,
$$
\zl=\frac{d \ul}{d\lm}=\ul^{'},\quad\wl=\frac{d \zl}{d\lm}=\zl^{'},\quad <f>=\ino \rl f,\quad \zlz=\zl-<\zl>.
$$
The derivative of $-g(\lm)=\lm<\zl>-1$ takes the form,
\beq\label{crux0}
-g^{'}(\lm)=(\lm <\zl>)^{'}=<\zl>+\lm <\zl>^{'}
\eeq
or else, since $<\zl>^{'}=<\zlz^2>+<\wl>$,
\beq\label{alg1}
-g^{'}(\lm)=<\zl>+\lm <\zlz^2>+\lm<\wl>.
\eeq

Clearly we have
\beq\label{zl}
-\Delta \zl = \rl +\lm\rl \zlz.
\eeq

Since $\prl$ is a constrained problem, $\sg_{1,\lm}>0$ does not imply that the maximum principle holds for the linearized problem \eqref{2.1} relative to $\prl$. Therefore, we can not claim that $\sg_{1,\lm}>0$ implies $\zl\geq0$. However, as a consequence of Proposition \ref{positive}$(i)$ we are able to prove the following version of the maximum principle based on the sign of $g(\lm)$:
\bpr\label{prcrux}$\,$\\
$(i)$ $<\zl>>0$, $\forall\,\lm\in (-\ii,8\pi)$.\\
$(ii)$ If $\lm\in(-\ii, 4\pi)$, then $g(\lm)=1-\lm<\zl>> 0$.\\
$(iii)$ If $\lm\in (-\ii,8\pi)$ and $1-\lm<\zl>\geq 0$, then $\zl\geq 0$.
\epr
\proof$\,$\\
$(i)$ Multiplying \rife{zl}  by $\zl$ we conclude that,
\beq\label{zlz}
\ino |\nabla \zl|^2-\lm\ino \rl\zlz^2=<\zl>,
\eeq
that is, since by Proposition \ref{positive}$(i)$ the first eigenvalue $\sg_1$ of $L_{\ssb}$ is strictly positive for $\lm\in (-\ii,8\pi)$, then we have,
$$
<\zl>\geq\sg_1<\zlz^2>>0,\;\forall\,\lm\in (-\ii,8\pi).
$$

$(ii)$ Multiplying \rife{zl}  by $\zl$ we also conclude that,
$$
0< \nu_1<\zl^2>\leq \ino |\nabla \zl|^2-\lm<\zl^2>=<\zl>(1-\lm<\zl>),
$$
where $\nu_1$ is the first eigenvalue of $-\Delta \phi-\lm\rl \phi$, which for $\lm< 4\pi$ satisfies $\nu_1> 0$, see \cite{BLin3} and \cite{suz}.
Therefore, since $<\zl>>0$, then if $\lm< 4\pi$ we conclude that $1-\lm<\zl>> 0$.\\

$(iii)$ We argue by contradiction and suppose that $z_{-}=\zl \chi_{\om_-}\not\equiv 0$, where $\om_-=\{x\in \om\,:\, \zl<0\}$.
Here $\chi_A$ is the characteristic function of the set $A$. Put $z_{+}=\zl \chi_{\om_+}$ and
$\phi=\al z_{+}+z_{-}$. Then $\phi$ satisfies
$$
-\Delta \phi=\al \lm\rl z_+ +\al\rl \chi_+ +\lm\rl z_- + \rl \chi_- -\al\lm\rl <\zl>\chi_+-\lm\rl<\zl>\chi_-=
$$
$$
\lm\rl \phi-\al\lm\rl<z_+>\chi_+-\al\lm\rl<z_->\chi_+-\lm\rl<z_->\chi_- -\lm\rl<z_+>\chi_-+\al\rl\chi_+ +\rl\chi_-
$$
Multiplying this equation by $\phi$ and integrating by parts we find that,
$$
\ino |\nabla \phi|^2=\lm\ino\rl \phi^2-\lm<\al z_+>^2-\lm< z_->^2-\lm(\al^2+1)<z_-><z_+>+\al^2 <z_+> +<z_->.
$$
Since
$$
\left(\ino\rl \phi\right)^2=<\al z_+>^2 +<z_->^2+2\al<z_-><z_+>,
$$
we obtain that
$$
\ino |\nabla \phi|^2-\lm<\phi_0^2>=P(\al):=-\lm(\al-1)^2<z_+><z_->+\al^2<z_+>+<z_->.
$$
Then, by using once more that $\sg_1$ is strictly positive, we conclude that the r.h.s of this inequality, which we call $P(\al)$, must be strictly positive if $\lm< 8\pi$. Indeed, we see that,
\beq\label{firsteig}
P(\al)=\ino |\nabla \phi|^2-\lm\ino \rl<\phi^2_0>\geq \sg_1<\phi_0^2>,
\eeq
and that, since $z_-$ does not vanish identically by assumption, $<\phi_0^2>>0$ for any $\al\in \R$.

Therefore we have,
$$
P(\al)=<z_+>(1-\lm<z_->)\al^2+2\lm<z_+><z_->\al+<z_->(1-\lm<z_+>)>0, \forall\,\al\in\R.
$$
The determinant of $P(\al)$ takes the form $-4<z_+><z_->(1-\lm<\zl>)$
which is non negative whenever  $(1-\lm<\zl>)\geq 0$. Therefore, if $(1-\lm<\zl>)\geq 0$, then $P(\al)$ must be
non positive somewhere, which is impossible for $\lm<8\pi$. Therefore $z_-\equiv 0$ whenever  $(1-\lm<\zl>)\geq 0$ and $\lm<8\pi$,
as claimed.
\finedim

\bigskip

It turns out that the analysis of the sign of $g(\lm)$ is far from being trivial. Surprisingly enough, we will succeed in carrying out the latter analysis  after proving that $g(\lm)$ satisfies a first order O.D.E.
\ble
The function $g(\lm)=1-\lm<\zl>$ satisfies
\beq\label{cruxg}
g^{'}(\lm)=a(\lm)g(\lm)+b(\lm), \quad\lm \in(-\ii,8\pi),\quad \quad g(0)=1,
\eeq
where
$$
a(\lm)=-(2\lm<\zlz^2>+\lm<\zl^2>+<\zl>), \quad b(\lm)=-\lm^2<\zl^3>.
$$
\ele
\proof
From \rife{zl} we find that,
\beq\label{wl}
-\Delta \wl=2\rl \zlz+\lm\rl(\zlz^2)_{\sscp 0}+\lm\rl \wlz.
\eeq
Multiplying this equation by $\zl$, integrating by parts and using \rife{zl} we also find that,
\beq\label{crux1}
<\wl>=2<\zlz^2>+\lm <(\zlz^2)_{\sscp 0}\zlz>\equiv 2<\zlz^2>+\lm <\zlz^3>
\eeq

By using \rife{crux1} in \rife{alg1}, we reduce \rife{alg1} to,
\beq\label{crux2}
-g^{'}(\lm)=<\zl>+3\lm <\zlz^2>+\lm^2 <\zlz^3>.
\eeq

Next, observing that
$$
<\zlz^3>=<\zl^3>-3<\zl^2><\zl>+2<\zl>^3,
$$
we see that \rife{crux2} takes the form
$$
-g^{'}(\lm)=\lm^2 <\zl^3>-3\lm^2<\zl^2><\zl>+2\lm^2<\zl>^3+3\lm<\zl^2>-3\lm<\zl>^2+<\zl>=
$$
$$
\lm^2 <\zl^3>+3\lm<\zl^2>(1-\lm<\zl>)+2\lm^2<\zl>^3-2\lm <\zl>^2-\lm <\zl>^2+<\zl>=
$$
$$
\lm^2 <\zl^3>+3\lm<\zl^2>(1-\lm<\zl>)-2\lm<\zl>^2(1-\lm<\zl>)+<\zl>(1-\lm<\zl>)=
$$
$$
\lm^2 <\zl^3>+2\lm<\zlz^2>(1-\lm<\zl>)+\lm<\zl^2>(1-\lm<\zl>)+<\zl>(1-\lm<\zl>).
$$
In particular we see that,
$$
-g^{'}(\lm)=\lm^2 <\zl^3>+2\lm<\zlz^2> g(\lm)+\lm<\zl^2>g(\lm)+<\zl>g(\lm),\quad g(0)=1,
$$
for any $\lm\in (-\ii,8\pi)$, as claimed. \finedim

\bigskip

At this point we can conclude the proof of Lemma \ref{lem1.3}.\\
By \rife{cruxg} we find that
\beq\label{cruxg1}
e^{-A(\lm)}g(\lm)=1+\int\limits_{0}^{\lm}e^{-A(t)}b(t)dt=1-\int\limits_{0}^{\lm}e^{-A(t)}t^2<z_{\sscp t}^3>dt,
\eeq
where $A(\lm)=\int_{0}^{\lm}a(t)dt$. By Proposition \ref{prcrux}$(ii)$ we have $g(\lm)>0$ for $\lm< 4\pi$. Since $\zl=\ul^{'}=(\lm\pl)^{'}$, then
$<\zl>=<\pl>+\lm<\etb>=2\el+\lm \el^{'}$, where we used \rife{elm} and \rife{8.12.2}. Thus, as a consequence of \rife{delm}, we also conclude that
 $$
 <\zl>\geq 2\el\to +\ii,\mbox{ as }\lm\to 8\pi^{-},
 $$
where we used \rife{elm1}. Therefore $g(\lm)\to -\ii$ as $\lm\to 8\pi^{-}$.
Since $g(0)=1$ and $g(\lm)$ is continuous, then there exists at least one value of $\lm=\lm_0\in[4\pi,8\pi)$ such that $g(\lm_0)=0$. Let
$$
\lm_*=\sup\{\lm>0\,:\,g(\tau)>0,\,\forall\,0\leq \tau<\lm\}.
$$
Then $\lm_*\geq 4\pi$, $g(\lm_*)=0$ and we claim that $g(\lm)<0$ for any $\lm>\lm_*$. Indeed, by Proposition~\ref{prcrux}~$(iii)$ we have $\zl\geq 0$ for  $\lm\leq \lm_*$. Therefore, by continuity, $<\zl^3>>0$ in a small enough right neighborhood of $\lm_*$ and so, by \rife{cruxg1},
$g(\lm)<0$ in  a small enough right neighborhood of $\lm_*$ as well. We argue by contradiction and suppose that the claim is false. Therefore,
putting
$$
\widehat{\lm}=\sup\{\lm>\lm_*\,:\,g(\tau)<0,\,\forall\,\lm_*<\tau<\lm\},
$$
we must have that $\widehat{\lm}\in(\lm_*,8\pi)$. Since clearly $g(\widehat{\lm})=0$ and $g(\lm)<0$ for $\lm_*<\lm<\widehat{\lm}$, then $g^{'}(\widehat{\lm})\geq 0$ and then by \rife{cruxg}
$$
(\widehat{\lm})^2\left.<\zl^3>\right|_{\lm=\widehat{\lm}}=-g^{'}(\widehat{\lm})\leq 0.
$$

However this is impossible since by Proposition \ref{prcrux} $(iii)$ we have $\left.\zl\right|_{\lm=\widehat{\lm}}\geq 0$ which implies
that  $\left.<\zl^3>\right|_{\lm=\widehat{\lm}}>0$. This contradiction shows that in fact $g(\lm)<0$ for $\lm>\lm_*$ as claimed. \finedim

\bigskip

\section{Appendix} \label{sec:app}

We present an example about the spectral analysis introduced in section \ref{sec:spectral}. To this end we consider a simplified linear problem which however share the same structure as \rife{2.1}:
\beq\label{2.1.0}
-\Delta \phi=\sg\left(\phi-\int\limits_{B_1}\hspace{-.4cm}-\phi\right) \quad\mbox{ in }\quad B_1
\eeq
with $\phi=0$ on $\partial B_1$, and $\int\limits_{B_1}\hspace{-.4cm}-\phi=\frac{1}{|B_1|}\int\limits_{B_1}\phi$. Passing to the new variable $\phi_0=\phi-\int\limits_{B_1}\hspace{-.4cm}-\phi$ we are reduced to
calculate the spectrum of
\begin{equation}\label{D.1}
-\Delta \phi_0-\sg \phi_0=0\quad\mbox{in}\quad B_1,
\end{equation}
with boundary conditions
\begin{equation}\label{D.2}
\mbox{(I)}\quad\phi_0=\mbox{constant}\quad\mbox{on}\quad \partial B_1\quad\mbox{and}\quad \mbox{(II)}\quad \int\limits_{B_1}\hspace{-.4cm}-\phi_0=0.
\end{equation}
Passing to polar coordinates, we first consider solutions of (\ref{D.1}) of the form
$$
\psi_{n}(r,\theta)=\left(A\cos{(n\theta)}+B\sin{(n\theta)}\right)J_{n}(\sqrt{\sigma}r)
$$
where $J_n$ is the Bessel function of order $n$.
Now either $n=0$ and then (\ref{D.2})-(I) is always satisfied or $n\geq 1$ and then
(\ref{D.2})-(I) is satisfied if and only if $\sigma=\sigma_{n,m}=\mu_{n,m}^2$, where $\mu_{n,m}$ is the m-th zero of
$J_n$. Next, if $n\geq 1$ then (\ref{D.2})-(II) is always satisfied, while if $n=0$ then (\ref{D.2})-(II) is equivalent to
$\int\limits_0^{1}J_0(\sqrt{\sigma}r)rdr=0$. Since $(J_1(r)r)^{'}=rJ_0(r)$ this is equivalent to
$\sqrt{\sigma}J_1(\sqrt{\sigma})=0$, that is $\sigma=\mu_{1,m}^2$. Therefore the eigenvalues of \rife{2.1.0} are
$$
\sigma_{n,m}=\mu_{n,m}^2,
$$
where $\sigma_{1,m}$ admits three eigenfunctions, $\{J_0(\mu_{1,m}r)-J_0(\mu_{1,m}),\cos{(\theta)}J_1(\mu_{1,m}r),\sin{(\theta)}J_1(\mu_{1,m}r)\}$  and  $\sigma_{n,m}$ with $n\geq 2$ that admits two eigenfunctions $\{\cos{(n\theta)}J_n(\mu_{n,m}r),\sin{(\theta)}J_n(\mu_{n,m}r)\})$.  Observe that these are
the eigenfunctions of \rife{2.1.0}. In particular the first eigenvalue
$\sigma_{1,1}=\mu_{1,1}^2\simeq (3.83)^2$ admits three eigenfunctions, one of which is radial,
$$
\phi_1=\phi_1(r)=J_0(\mu_{1,1}r)-J_0(\mu_{1,1}),
$$

and satisfies $\phi_1(1)=0$ and $\phi_1^{'}(1)=\mu_{1,1}J_0^{'}(\mu_{1,1})=0$. This is not in contradiction with the Hopf Lemma or with general unique continuation principles, since  $\int\limits_{B_1}\hspace{-.4cm}-\phi_1=-J_0(\mu_{1,1})\neq 0$ and the identically zero function $\phi\equiv 0$ is not a solution of \rife{2.1.0} with $\int\limits_{B_1}\hspace{-.4cm}-\phi\neq 0$. Alternatively, in terms of $\phi_0$, the eigenfunction $\phi_{1,0}=\phi_1-\int\limits_{B_1}\hspace{-.4cm}-\phi_1$ satisfies $\phi_{1,0}(1)=J_0(\mu_{1,1})\neq 0$ and $\phi^{'}_{1,0}(1)=\mu_{1,1}J_0^{'}(\mu_{1,1})=0$ but still the function $\phi_0\equiv J_0(\mu_{1,1})\neq 0$ is not a solution of \rife{D.1}, \rife{D.2}.


\begin{thebibliography}{99}

\bibitem{Ban1}
C. Bandle, {\em Existence theorems, qualitative results and a priori bounds for a class of nonlinear Dirichlet problems},
Arch. Rational Mech. Anal. {\bf 58}(3) (1975),  219-238.

\bibitem{bav} F. Bavaud, {\em Equilibrium properties of the Vlasov functional: the generalized Poisson-Boltzmann-Emden
equation}, Rev. Mod. Phys. {\bf 63}(1) (1991), 129-149.


\bibitem{B2} D. Bartolucci, {\em Stable and unstable equilibria of
uniformly rotating self-gravitating cylinders}, {Int. Jour. Mod. Phys. D} {\bf 21}(13) (2012), 1250087.


\bibitem{Bons} D. Bartolucci, {\em Global bifurcation analysis of mean field equations and
the Onsager microcanonical description of two-dimensional turbulence}, {Calc. Var. \& P.D.E.} {\bf 58}:18 (2019).

\bibitem{BdM2} D. Bartolucci, F. De Marchis,
{\em Supercritical Mean Field Equations on convex domains and the Onsager's
statistical description of two-dimensional turbulence}, Archive for Rational Mechanics and Analysis {\bf 217}/2 (2015), 525-570;
DOI: 10.1007/s00205-014-0836-8.

\bibitem{BJLY3} D. Bartolucci, A. Jevnikar, Y. Lee, W. Yang,
{\em Local uniqueness of m-bubbling sequences for the Gel'fand equation,} Comm. in P.D.E. \textbf{44}(6) (2019), 447-466.

\bibitem{BLin3} D. Bartolucci, C.S. Lin, {\em Existence and uniqueness for
Mean Field Equations on multiply connected domains at the critical parameter},
{Math. Ann.} {\bf 359} (2014), 1-44; DOI 10.1007/s00208-013-0990-6.

\bibitem{BMal} D. Bartolucci, A. Malchiodi, {\em Mean field equations and domains of first kind}, Preprint 2020,
arXiv:2001.01654.

\bibitem{Bw}  D. Bartolucci, G.Wolansky, {\em Maximal entropy solutions under prescribed mass and energy}, 
Jour. Diff. Eq., {\bf 268} (2020) 6646-6665; DOI: 10.1016/j.jde.2019.11.040.

\bibitem{beb} J. Bebernes, D. Eberly, Mathematical Problems from Combustion Theory, A. M.
S. 83, Springer-Verlag, New York (1989).

\bibitem{BKN} P. Biler, A. Krzywicki, T. Nadzieja, {\em Self-interaction of brownian particles
coupled with thermodynamic processes}, Rep. Math. Phys. {\bf 42} (1998), 359-372.

\bibitem{bm} H. Brezis, F. Merle,
{\em Uniform estimates and blow-up behaviour for solutions of $-\Delta u = V(x)e^{u}$ in two dimensions},
Comm. in P.D.E. {\bf 16}(8,9) (1991), 1223-1253.

\bibitem{bdt} B. Buffoni, E.N. Dancer, J.F. Toland, {\em The sub-harmonic bifurcation of Stokes waves},
Arch. Rat. Mech. Anal. {\bf 152}(3) (2000), 24-271.

\bibitem{but} B. Buffoni, J. Toland, {Analytic Theory of Global Bifurcation}, (2003) Princeton Univ. Press.

\bibitem{clmp2} E. Caglioti, P.L. Lions, C. Marchioro, M. Pulvirenti,
{\em A special class of stationary flows for two dimensional Euler equations: a
statistical mechanics description. II}, Comm. Math. Phys. {\bf 174} (1995),
229-260.

\bibitem{CCL} S.Y.A. Chang, C.C. Chen, C.S. Lin, {\em Extremal functions for a mean field equation in two dimension},
in: "Lecture on Partial Differential Equations", New Stud. Adv. Math. {\bf 2} Int. Press, Somerville, MA, 2003, 61-93.


\bibitem{CLin1} C. C. Chen, C.S. Lin, {\em Sharp estimates for solutions of multi-bubbles in compact Riemann surfaces},
                Comm. Pure Appl. Math. {\bf55} (2002), 728-771.

\bibitem{CLin2} C. C. Chen, C.S. Lin, {\em Topological Degree for a mean field equation on Riemann surface},
                Comm. Pure Appl. Math. {\bf 56} (2003), 1667-1727.



\bibitem{CrRab} M. G. Crandall, P. H. Rabinowitz, {\em Some Continuation and Variational
Methods for Positive Solutions of Nonlinear Elliptic Eigenvalue Problems},
{Arch. Rat. Mech. An.} {\bf  58} (1975), 207-218.



\bibitem{Dan} N. Dancer, {\em Global structure of the solutions of non-linear real analytic eigenvalue problems},
Proc. Lond. Math. Soc. (3) {\bf 27} (1973) , 747-765.

\bibitem{Dan2} N. Dancer, {\em Finite Morse index solutions of exponential problems},
Ann. I. H. P. An. Nonlin. (3) {\bf 25} (2008), 173-179.

\bibitem{Dan3} N. Dancer {\em Finite Morse index solutions of supercritical problems},
 Jour Reine Angew. Math. {\bf 620} (2008) , 213-233.

\bibitem{DanF} N. Dancer, A Farina {\em On the classification of solutions of $-\Delta u =e^{u}$
on $\R^n$ : stability outside a compact set and applications},
 Proc. A.M.S.. {\bf 137} (2009) , 1333-1338.

\bibitem{Gel} I. M. Gelfand, {\em Some problems in the theory of quasi-linear equations},
Amer. Math. Soc. Transl. {\bf 29}(2) (1963), 295-381.

\bibitem{gnn} B. Gidas, W.M. Ni, L. Nirenberg, {\em Symmetry and related
properties via the maximum principle}, Comm. Math. Phys. {\bf 68} (1979),
209-243.

\bibitem{GOS} M. Grossi, H. Ohtsuka, T. Suzuki, {\em Asymptotic non-degeneracy of the multiple blow-up solutions of the Gel'fand problem in two space dimensions}, Adv. Diff. Eqs. {\bf 16}(1-2) (2011), 145-164.

\bibitem{JL} D.D. Joseph, T.S. Lundgren, {\em
Quasilinear Dirichlet problems driven by positive sources}, Arch. Ration. Mech. Anal. {\bf 49} (1972/73) 241-269.


\bibitem{HK} M. Holzmann, H. Kielh\"ofer, {\em Uniqueness of global positive solution branches
of nonlinear elliptic problems,} Math. Ann. {\bf  300}, 221-241 (1994).



\bibitem{KLB}
J. Katz, D Lynden-Bell, \emph{The Gravothermal instability in two dimensions,}  M.N.R.A.S. \textbf{184} (1978) 709-712.

\bibitem{KK}
 J.P. Keener,  H.B. Keller,  \emph{Positive solutions of convex nonlinear eigenvalue problems,}
 J. Differential Equations \textbf{16} (1974), 103-125.

\bibitem{KC} H.B. Keller, D.S. Cohen, {\em Some positone problems suggested by nonlinear heat generation}, Jour. Math. Mech. {\bf 16} (1967).

\bibitem{KMdP} M. Kowalczyk, M. Musso, M. del Pino,  {\em Singular limits in
Liouville-type equations}, Calc. Var. {\bf 24}(1)  (2005), 47-81.

\bibitem{Kor}  P. Korman, {\sl Global Solution Curves for Semilinear Elliptic Equations}, World Scientific (2012).

\bibitem{Kor2} P. Korman, {\em Global solution curves for self-similar equations}, J. Differential Equations \textbf{257} (2014) 2543-2564.


\bibitem{KLO} P. Korman, Y. Li b, T. Ouyang, {\em A simplified proof of a conjecture for the perturbed
Gelfand equation from combustion theory}, J. Differential Equations \textbf{263} (2017) 2874-2885.




\bibitem{yy} Y.Y. Li,  {\em Harnack type inequality: the method of moving planes},
Comm. Math. Phys.  {\bf 200} (1999), 421--444.

\bibitem{Lions} P.L. Lions, \emph{On the Existence of Positive Solutions of Semilinear Elliptic Equations,}
S.I.A.M. Review \textbf{24}(4) (1982) 441-467.





\bibitem{moser} J. Moser, {\em A sharp form of an inequality by N.Trudinger}, Indiana Univ. Math. J. 20 (1971), 1077-1091.


\bibitem{NS90} K. Nagasaki, T. Suzuki, {\em Asymptotic analysis for two-dimensional elliptic eigen-
values problems with exponentially dominated nonlinearities}, Asymptotic Analysis,
{\bf 3} (1990), 173-188.

\bibitem{Os} J. Ostriker, \emph{The equilibrium of politropic and isothermal cylinders}, Astrophys. J. \textbf{140} (1964) 1056-1066.


\bibitem{Rab} P.H. Rabinowitz, {\em Some global results for nonlinear eigenvalue problems}, J. Funct. Anal.
{\bf 7} (1981), 487-513.

\bibitem{suz} T. Suzuki, {\em Global analysis for a two-dimensional elliptic eiqenvalue problem with the exponential
                nonlinearly}, Ann. Inst. H. Poincar\'e Anal. Non Lin\'eaire {\bf 9}(4) (1992), 367-398.

\bibitem{suzC} T. Suzuki, "Free Energy and Self-Interacting Particles", PNLDE
{\bf 62}, Birkhauser, Boston, (2005).





\end{thebibliography}
\end{document}